\documentclass[reqno,12pt]{amsart}
\usepackage{indentfirst}
\usepackage{amsmath}
\usepackage[T1]{fontenc}
\usepackage{amssymb}			
\usepackage{amsthm}
\usepackage{amsfonts}
\usepackage{enumerate}
\usepackage{enumitem}
\usepackage{lmodern}
\usepackage{setspace}
\usepackage{mathtools}
\usepackage[body={7in,9.5in},top=0.5in, left=0.8in ]{geometry}
\newcommand{\bn}{\mathbb{N}}

\newcommand{\br}{\mathbb{R}}
\newcommand{\ov}[1]{\overline{#1}}

\newcommand{\pref}[1]{\textup{(\ref{#1})}}

\newtheorem{thrm}{Theorem}%[section] %definicje wlasne; kolejno 
\newtheorem{df}[thrm]{Definition}		%numerowane w ramach sekcji
\newtheorem{lm}[thrm]{Lemma}

\newtheorem{rmk}[thrm]{Remark}

\newtheorem{xmpl}[thrm]{Example}

\author[D.Bugajewski]{Dariusz Bugajewski}
\address{D. Bugajewski, Department of Nonlinear Analysis and Applied Topology\\
  Faculty of Mathematics and Computer Science\\
  Adam Mickiewicz University in Pozna\'n\\
  ul.\ Umultowska 87\\
  61-614 Pozna\'n\\
  Poland}
\email{ddbb@amu.edu.pl}

\author[P.Ma\'ckowiak]{Piotr Ma\'ckowiak} 
\address{P.Ma\'ckowiak, Department of Nonlinear Analysis and Applied Topology\\
  Faculty of Mathematics and Computer Science\\
  Adam Mickiewicz University in Pozna\'n\\
  ul.\ Umultowska 87\\
  61-614 Pozna\'n\\
  Poland}
\email{piotr.mackowiak@amu.edu.pl}

\date{}
\title{A fixed point theorem for mappings in partial metric spaces}
\begin{document}
\begin{abstract}
In this paper we are going to prove a very general fixed point theorem for mappings acting in partial metric spaces. In that theorem we impose some conditions on behavior of considered mappings on orbits and a condition relating orbits of points of small size.
\end{abstract}
\maketitle 

\section{Introduction} 

In the paper we consider some class of mappings acting in partial metric spaces. Let us recall that the notion of a partial metric space was introduced by Matthews \cite{M} and it has been applied, among others, in the study of denotational semantics of a programming language. Other examples of partial metrics are connected with weighted graphs (see \cite{ST}, Section 3 for more details) or with fixed point theory (see \cite{BKMP}, \cite{IPR} and \cite{M}). 

Obviously, partial metric spaces are proper generalization of metric spaces, however it seems to be worth noticing that these classes of spaces possess some common properties. For example, as it was shown in \cite{BMW}, in the case of partial metric spaces the notion of compactness is equivalent to sequential compactness. 

The fixed point theory in partial metric spaces seems to be well developed. The partial metric version of the Banach Contraction Principle was proved in \cite{M}. Various extensions and generalizations of that principle were established in numerous papers (see e.g. \cite{ASS}, \cite{BW}, \cite{CKT}, \cite{IPR}, \cite{IPR2}, \cite{R} and \cite{ZZ}). In particular, the extension given in \cite{IPR} turned out to be very interesting as it enabled to associate completeness of partial metric spaces with the continuum hypothesis (see \cite{Mack2023} for details). 

 Let us emphasize that there are several partial metric fixed point theorems that can be obtained from the corresponding results valid in metric spaces, as it was noticed in the article \cite{HRS}.  

The goal of our article is to prove a very general fixed point theorem for mappings acting in partial metric spaces. According to our knowledge that new result has no preceding metric analogue. The main assumptions of our fixed point theorem describe the behavior of iterations of a mapping under consideration in a special way. 

Our note is organized as follows. In Section 2 we collect necessary definitions and notations. Next, in Section 3, we present some lemmas and the main result. Finally, in Section 4, we illustrate that new fixed point theorem by a few examples and compare it with some results known from the literature. In particular, we show that our result is a generalization of the Generalized Banach Contraction Principle proved in the paper \cite{IPR}. We also give examples showing that none of our main conditions ensures the existence of a fixed point without satisfying the second condition.

\section{Preliminaries} 

We denote by $\bn$ the set of positive integers, $\bn_0:=\bn\cup\{0\}$. $\br$ stands for the set of real numbers while $\br_+$ is the set of nonnegative real numbers. For any function $T:A\to A$ let $T^0(x):=x,\, x\in A$, where $A\neq\emptyset$. 

Because we are going to consider mappings acting in partial metric spaces, in what follows we recall some basic definitions and facts describing this class of spaces. 
\begin{df}
Let $U$ be a nonempty set. The mapping $p:U\times U\rightarrow \mathbb{R}_+$ is said to be a \emph{partial metric} or \emph{pmetric}  on $U$ if, for any $x,y,z\in U$, the following four conditions hold:

{\rm (P1)} $x=y$ if and only if $p(x,x)=p(y,y)=p(x,y)$;

{\rm (P2)} $p(x,x)\leq p(x,y)$;

{\rm (P3)} $ p(x,y)=p(y,x)$;

{\rm (P4)} $p(x,y)\leq p(x,z)+p(z,y)-p(z,z)$.

The set $U$ with the partial metric $p$ is called a \emph{partial metric space}, which is denoted by $(U,p)$. We call $p(x,x)$ the \emph{size} of $x$. 
\end{df} 
Obviously, a metric space is also a partial metric space with all points of size $0$.

Let $(U,p)$ be a partial metric space. A sequence $(x_q)_{q\in \mathbb{N}}$ of elements of the partial metric space $U$ converges to an element $x\in U$ if and only if
\[
\lim_{q\rightarrow \infty}p(x_q,x)=p(x,x),
\]
and we say that $x$ is a limit point of the sequence $(x_q)_{q\in \mathbb{N}}$. Let us notice that the limit point of the sequence $(x_q)_{q\in \mathbb{N}}$ in partial metric space need not be unique. 

A sequence $(x_q)_{q\in \mathbb{N}}$ of elements of a partial metric space $(U,p)$ is said to be a \emph{Cauchy sequence} if there exists the limit 
\[
\lim_{q,{q'}\rightarrow\infty}p(x_q,x_{q'}).
\]
In particular, a sequence $(x_q)_{q\in \mathbb{N}}$ of elements of a partial metric space $(U,p)$ is said to be an $a-$\emph{Cauchy  sequence} ($a\in \mathbb{R}_+)$ if
\[
\lim_{q,{q'}\rightarrow\infty}p(x_q,x_{q'})=a.
\]
A partial metric space $(U,p)$ is said to be \emph{complete} if every Cauchy sequence $(x_q)_{q\in \mathbb{N}}$ of elements of $U$ converges to an element $x\in U$ such that
\[
\lim_{q,{q'}\rightarrow\infty}p(x_q,x_{q'})=\lim_{q\rightarrow\infty}p(x_q,x)=p(x,x).
\]
Analogously, a partial metric space $(U,p)$ is said to be $a-$\emph{complete} if every $a-$Cauchy sequence $(x_q)_{q\in \mathbb{N}}$ of elements of $U$ converges to an element $x\in U$ such that
\[
\lim_{n\rightarrow\infty}p(x_q,x)=p(x,x)=a.
\]
As in \cite{IPR}, let us define $$\rho_p:=\inf\{p(x,x):\, x\in U\},\quad U_p:=\{x\in U:\, p(x,x)=\rho_p\}.$$
\section{Main results} 
Let $(U,p)$ be a partial metric space and let $\alpha\in [0,1)$ be fixed. In this section we are going to consider mappings $T:U\to U$ satisfying the following conditions:
\begin{enumerate}[label={\textup{(\Alph*)}},ref=\textup{\Alph*}]
\item\label{cnd:1:1} for $x\in U$, $q\in \bn$, it holds $$p(T(x),T^{q+1}(x))\leq \max\{\alpha p(x,T^{q+1}(x)),\ldots,\alpha p(x,T^2(x)),\alpha p(x,T(x)),p(x,x)\},$$ 
%\item\label{cnd:1:2} for any $\varepsilon>0,\, x\in U$, there exists $n\in \bn$ such that $$p(T^{n}(x),T^{n}(x))\leq p(x,x)+\varepsilon.$$
\item\label{cnd:1:2} there exists $\varepsilon_1>0$ such that, for any $x,y\in U$ with $\max\{p(x,x),p(y,y)\}\leq \rho_p+\varepsilon_1$, for each $\varepsilon>0$, there exists $q_\varepsilon\in \bn$ for which $$p(T^{q_\varepsilon}(x),T^{q_\varepsilon}(y))\leq \max\{p(x,x),p(y,y)\}+\varepsilon.$$
\end{enumerate}

Notice that both conditions \pref{cnd:1:1} and \pref{cnd:1:2} are satisfied if $T:U\to U$ is a contraction and $(U,p)$ is a metric space. It is also clear that that the inequality appearing in condition \pref{cnd:1:1} is satisfied by any fixed point of $T$ for any $\alpha\in [0,1)$. Moreover, condition \pref{cnd:1:2} is satisfied by any mapping for which there exists a unique fixed point $\ov{x}$ in the set $U_p$ and the sequence  $(T^q(x))_{q\in \bn}$ converges to $\ov{x}$ for every $x\in U$. We postpone further discussion of conditions \pref{cnd:1:1} and \pref{cnd:1:2} to Section \ref{sec:discussion}.

In what follows we will use some handy notation, namely, for any function $T:U\to U$ and point $x\in U$, we denote $$x_q:=T^q(x),\,q\in \bn_0.$$ Condition \pref{cnd:1:1} implies the following three crucial lemmas.

\begin{lm}\label{lm:1a} If a mapping $T:U\to U$ satisfies condition \pref{cnd:1:1} for some $\alpha\in [0,1)$, then $0\leq p(x_q,x_{q'})\leq \frac{2}{1-\alpha}p(x_0,x_1)$, $x\in U$, $q,q'\in \bn_0$.
\end{lm} 

\begin{proof} 
By \pref{cnd:1:1} and properties of the partial metric $p$, for $k,l\in \bn$, we get 
\begin{multline*}
p(x_1,x_{k+1})\leq \max\{\alpha p(x_0,x_{k+1}),\alpha p(x_0,x_{k}),\ldots,\alpha p(x_{0},x_{2}),\underbrace{\alpha p(x_{0},x_{1}),p(x_{0},x_{0})}_{\leq p(x_0,x_1)}\}\\\leq \max\{\alpha(p(x_0,x_1)+p(x_1,x_{k+1})), \alpha (p(x_0,x_1)+p(x_1,x_{k})),\ldots,\alpha(p(x_0,x_1)+p(x_1,x_2)),p(x_{0},x_{1})\}\\\leq p(x_0,x_1)+\max\{\alpha p(x_1,x_{k+1}), \alpha p(x_1,x_{k}),\ldots,\alpha p(x_1,x_2)\}\\=p(x_0,x_1)+\alpha \max\{p(x_1,x_{k+1}), p(x_1,x_{k}),\ldots,p(x_1,x_2)\}\\\leq p(x_0,x_1)+\alpha \max\{\alpha p(x_0,x_{k+1}), \alpha p(x_0,x_{k}),\ldots,\alpha p(x_0,x_1), p(x_0,x_0)\}\\\leq p(x_0,x_1)(1+\alpha)+\alpha^2 \max\{p(x_1,x_{k+1}),p(x_1,x_{k}),\ldots,p(x_1,x_2)\}\\\leq p(x_0,x_1)(1+\alpha+\alpha^2)+\alpha^3 \max\{p(x_1,x_{k+1}),p(x_1,x_{k}),\ldots,p(x_1,x_2)\}\leq\ldots\\\leq p(x_0,x_1)(1+\alpha+\ldots+\alpha^{l-1})+\alpha^l\max\{p(x_1,x_{k+1}),\ldots, p(x_1,x_2)\},
\end{multline*}
and, as $l\to+\infty$, $p(x_1,x_{k+1})\leq \frac{1}{1-\alpha}p(x_0,x_1)$. Thus $p(x_q,x_{q+k})\leq p(x_q,x_1)+p(x_1,x_{q+k})\leq \frac{2}{1-\alpha}p(x_0,x_1)$ for $q,k\in \bn_0$.
\end{proof}

\begin{lm}\label{lm:2a} If a mapping $T:U\to U$ satisfies condition \pref{cnd:1:1} for some $\alpha\in [0,1)$, then, for any $x\in U$, there exists the limit $\lim_{q\to +\infty} p(x_q,x_q)$ and the sequence $(x_q)_{q\in\bn_0}$ is $r_x$-Cauchy, where $r_x:=\lim_{q\to\infty}p(x_q,x_q)$.
\end{lm} 

\begin{proof}
For $q\in \bn,k\in\bn_0$, we get 
\begin{multline}\label{eq:0a}
p(x_{q},x_{k+q})\leq \max\{\alpha p(x_{q-1},x_{k+q}),\alpha p(x_{q-1},x_{k+q-1}),\ldots, \alpha p(x_{q-1},x_q), p(x_{q-1},x_{q-1})\}\\\leq \max\{\alpha^2 p(x_{q-2},x_{k+q}),\alpha^2 p(x_{q-2},x_{k+q-1}),\ldots, \alpha^2p(x_{q-2},x_{q-1}), \alpha p(x_{q-2},x_{q-2}),\, p(x_{q-1},x_{q-1})\}\\\leq\max\{\alpha^3 p(x_{q-3},x_{k+q}),\ldots, \alpha^3p(x_{q-3},x_{q-2}), \alpha^2 p(x_{q-3},x_{q-3}), \alpha p(x_{q-2},x_{q-2}),\, p(x_{q-1},x_{q-1})\}\leq \ldots\\\leq
\max\{\alpha^{q-1}p(x_{1},x_{k+q}),\ldots,\alpha^{q-1}p(x_{1},x_{2}),\alpha^{q-2} p(x_{1},x_{1}),\ldots,\alpha p(x_{q-2},x_{q-2}),p(x_{q-1},x_{q-1})\}.
\end{multline}
Let $r:=\limsup_{q\to \infty}p(x_q,x_{q+1})$, $r':=\liminf_{q\to \infty}p(x_q,x_{q})$. Obviously, $0\leq r'\leq r$ and, by Lemma~\ref{lm:1a}, $r<+\infty$. Fix now any $\varepsilon>0$ and choose $\ov{q}\in \bn$ such that $\alpha^{\ov{q}-1}\frac{2}{1-\alpha}p(x_0,x_1)\leq \alpha(r+\varepsilon)$. Let also $\ov{q}_0\in \bn$, $\ov{q}_0\geq 2\ov{q}$, be such that $r'-\varepsilon\leq p(x_{q-1},x_{q-1})\leq p(x_{q-1},x_{q})\leq r+\varepsilon$ for $q\geq \ov{q}_0$. Hence, by inequality (\ref{eq:0a}) and Lemma \ref{lm:1a}, for $q\geq 2\ov{q}_0$, $k\in \bn_0$, we obtain 
\begin{multline}\label{eq:0b} 
p(x_{q},x_{k+q})\leq \max\{\alpha^{q-1}p(x_{1},x_{k+q}),\ldots,\alpha^{{q}-1}p(x_{1},x_{2}),\alpha^{q-2} p(x_{1},x_{1}),\ldots,\alpha p(x_{q-2},x_{q-2}),p(x_{q-1},x_{q-1})\}\\=\max\{\alpha^{q-1}p(x_{1},x_{k+q}),\ldots,\alpha^{{q}-1}p(x_{1},x_{2}),\alpha^{q-2} p(x_{1},x_{1}),\ldots,\alpha^{\ov{q}-1} p(x_{q-\ov{q}},x_{q-\ov{q}}),\\\alpha^{\ov{q}-2} p(x_{q-\ov{q}+1},x_{q-\ov{q}+1}),\ldots,\alpha p(x_{q-2},x_{q-2}),p(x_{q-1},x_{q-1})\}\\
\leq \max\{\alpha^{\ov q-1}\frac{2}{1-\alpha}p(x_0,x_1),\alpha^{\ov{q}-2} p(x_{q-\ov{q}+1},x_{q-\ov{q}+1}),\ldots,\alpha p(x_{q-2},x_{q-2}),p(x_{q-1},x_{q-1})\}\\\leq \max\{\alpha(r+\varepsilon),\alpha^{\ov{q}-2} p(x_{q-\ov{q}+1},x_{q-\ov{q}+1}),\ldots,\alpha p(x_{q-2},x_{q-2}),p(x_{q-1},x_{q-1})\}\\\leq \max\{\alpha(r+\varepsilon),\alpha p(x_{q-\ov{q}+1},x_{q-\ov{q}+1}),\ldots,\alpha p(x_{q-2},x_{q-2}),p(x_{q-1},x_{q-1})\} \\\leq \max\{\alpha(r+\varepsilon),p(x_{q-1},x_{q-1})\}
\end{multline}
where we used the facts that $\alpha\in [0,1)$ and if $q\geq 2q_0$, then $q-\ov{q}\geq 2\ov q_0-\ov{q}\geq \ov q_0+2\ov q-\ov q\geq \ov{q}_0$.
Suppose that $r>r'\geq 0$. Let us now assume not only that $\varepsilon>0$, but $\varepsilon<\frac{1-\alpha}{1+\alpha}r$; observe that $\alpha(r+\varepsilon)<r-\varepsilon$. By the definition of $r$, there exists a subsequence $(q_k)_{k\in \bn_0}$ of $(n)_{n\in \bn_0}$ for which $p(x_{q_k},x_{q_k+1})\geq r-\varepsilon,\, k\in \bn_0$. Without loss of generality we may assume that $q_0\geq 2\ov{q}_0$. By inequality (\ref{eq:0b}) and the choice of $\varepsilon$, for $k\in \bn$, we have
$$r-\varepsilon\leq p(x_{q_k},x_{q_k+1})\leq \max\{\alpha (r+\varepsilon),p(x_{q_k-1},x_{q_k-1})\}=p(x_{q_k-1},x_{q_k-1}).$$
Thus 
\begin{multline*}r-\varepsilon\leq p(x_{q_k-1},x_{q_k-1})\leq p(x_{q_k-1},x_{q_k})\leq \max\{\alpha (r+\varepsilon),p(x_{q_k-2},x_{q_k-2})\}=p(x_{q_k-2},x_{q_k-2})\\\leq p(x_{q_k-2},x_{q_k-1})\leq \max\{\alpha (r+\varepsilon),p(x_{q_k-3},x_{q_k-3})\}=p(x_{q_k-3},x_{q_k-3})\leq\ldots\\\leq p(x_{q_{k-1}},x_{q_{k-1}})\leq p(x_{q_{k-1}},x_{q_{k-1}+1})\leq \max\{\alpha (r+\varepsilon),p(x_{q_{k-1}-1},x_{q_{k-1}-1})\}\leq r+\varepsilon,\end{multline*}
and, subsequently, for $k\in \bn$ and $q\in\{q_{k-1},q_{k-1}+1,q_{k-1}+2,\ldots,q_k\}$, 
$$r-\varepsilon\leq p(x_{q},x_{q})\leq p(x_{q},x_{q+1})\leq r+\varepsilon,$$
which, due to the fact that $\varepsilon>0$ can be arbitrarily close to $0$, implies that the limit $\lim_{q\to\infty}p(x_q,x_{q+1})$ exists and that inequality $r'<r$ is impossible. Therefore, $r=r'=\lim_{q\to \infty}p(x_q,x_q)$. Now, since formula (\ref{eq:0b}) with $r$ replaced with $r_x:=\lim_{q\to \infty}p(x_q,x_q)$ is true for any fixed $\varepsilon>0$ and all $k\in \bn$ and sufficiently large $q$, we obtain
$$r_x-\varepsilon\leq p(x_q,x_q)\leq p(x_{q},x_{k+q})\leq \max\{\alpha(r_x+\varepsilon),p(x_{q-1},x_{q-1})\}\leq r_x+\varepsilon,$$
for sufficiently large $q$ and all $k\in \bn$, which proves that the sequence $(x_q)_{q\in \bn_0}$ is $r_x$-Cauchy.
\end{proof}

\begin{lm}\label{lm:2} 
Suppose that a mapping $T:U\to U$ satisfies condition \pref{cnd:1:1} for some $\alpha\in [0,1)$. If $x\in U_p$ and $\lim_{q\to\infty}p(T^q(x),T^q(x))=\lim_{n\to\infty}p(T^q(x),x)=p(x,x)$, then $x=T(x)$.
\end{lm}
\begin{proof}
Assume, for some $x\in U_p$, we have $\lim_{q\to\infty}p(T^q(x),T^q(x))=\lim_{q\to\infty}p(T^q(x),x)=p(x,x)=\rho_p$. It suffices to show that $p(x,T(x))=\rho_p$. Suppose to the contrary, that is, we have $p(x,T(x))>\rho_p\geq 0$. Since $\lim_{n\to\infty}p(T^q(x),x)=\rho_p$, we get that $p(x,T^{\ov q}(x))=\max_{q\in \bn}p(x,T^q(x))\geq p(x,T(x))$ for some ${\ov q}\in \bn$ and may assume that $\ov{q}\geq 1$ is the smallest among such numbers; notice that $p(x,T^{\ov q-1}(x))<p(x,T^{\ov{q}}(x))$. Let $x_q:=T^q(x),\, q\in \bn$. 

As in the proof of Lemma \ref{lm:2a} (see formula \pref{eq:0a}) we get, for $q,k\in \bn$,
$$p(x_{q},x_{k+q})\leq \max\{\alpha^{q-1}p(x_{1},x_{k+q}),\ldots,\alpha^{q-1}p(x_{1},x_{2}), \alpha^{q-2} p(x_{1},x_{1}),\ldots,\alpha p(x_{q-2},x_{q-2}),p(x_{q-1},x_{q-1})\}.$$
According to the above formula, by the definition of ${\ov q}$ and condition \pref{cnd:1:1}, we have, for all $k\in \bn$,
\begin{multline*}p(x_0,x_{\ov q})\leq \underbrace{p(x_0,x_{k+\ov q})-p(x_{k+\ov q},x_{k+\ov q})}_{p_{k}^{\ov q}:=}+p(x_{{\ov q}},x_{k+\ov q})
\\\leq p^{\ov q}_k+\max\{\alpha^{{\ov q}-1}p(x_{1},x_{k+{\ov q}}),\ldots,\alpha^{{\ov q}-1}p(x_{1},x_{2}),\alpha^{{\ov q}-2} p(x_{1},x_{1}),\ldots,\alpha p(x_{{\ov q}-2},x_{{\ov q}-2}),p(x_{{\ov q}-1},x_{{\ov q}-1})\}\\\leq 
p^{\ov q}_k+\max\{\alpha^{{\ov q}}p(x_{0},x_{k+{\ov q}}),\ldots, \alpha^{{\ov q}}p(x_{0},x_{1}),\alpha^{{\ov q}-1} p(x_0,x_0),\alpha p(x_{0},x_{{\ov q}}), p(x_{0},x_{{\ov q-1}})\}\\
\leq 
p^{\ov q}_k+\max\{\alpha p(x_{0},x_{{\ov q}}), p(x_{0},x_{{\ov q-1}})\}=(\star).
\end{multline*}
Now, by $p^{\ov q}_k=p(x,T^{k+\ov{q}}(x))-p(T^{k+\ov{q}}(x),T^{k+\ov{q}}(x))\to 0$, as $k\to \infty$, we have  
$$p(x_0,x_{\ov q})\leq\max\{\alpha p(x_{0},x_{{\ov q}}), p(x_{0},x_{{\ov q-1}})\},$$
which implies either $p(x_0,x_{\ov q})=0\leq \rho_p$ or $p(x_0,x_{\ov q})\leq p(x_{0},x_{{\ov q-1}})$, but the former is impossible by assumption, while the latter cannot hold due to the choice of $\ov{q}$. The claim follows.
\end{proof} 

The next lemma uses both conditions \pref{cnd:1:1} and \pref{cnd:1:2}.
\begin{lm}\label{lm:1} Suppose that a mapping $T:U\to U$ satisfies condition \pref{cnd:1:1} for some $\alpha\in [0,1)$ and condition \pref{cnd:1:2} for some $\varepsilon_1>0$. For any $x\in U$ such that $p(x,x)<\rho_p+\varepsilon_1$, we have  $r_x\leq p(x,x)$ and the sequence $(T^q(x))_{q\in\bn_0}$ is $r_x$-Cauchy, where $r_x:=\lim_{q\to\infty}p(T^q(x),T^{q}(x))$.
\end{lm} 

\begin{proof} 
Fix $x\in U$ with $p(x,x)<\rho_p+\varepsilon_1$. By Lemma \ref{lm:2a}, the sequence $(x_q)_{q\in \bn_0}$, where $x_q:=T^q(x),\, q\in \bn_0$, is $r_x$-Cauchy, with $r_x=\lim_{q\to\infty}p(T^q(x),T^{q}(x))$.
We shall now show that $r_x\leq p(x,x)$. To this end it suffices to show that, for each $\varepsilon\in (0,\, \rho_p+\varepsilon_1-p(x,x))$, there exists %$\ov{q}_\varepsilon\in \bn$ with $r_x-\varepsilon\leq p(x_q,x_q)$, $q\geq \ov{q}_\varepsilon$; notice that we may assume that 
arbitrarily large $q_\varepsilon\in \bn$ for which $p(x_{q_\varepsilon},x_{q_\varepsilon})<p(x,x)+\varepsilon$. By condition \pref{cnd:1:2}, there exists a nonempty set $M$ of all $q\in \bn$ for which $p(x_{q},x_{q})<p(x,x)+\varepsilon$. Assume that the set $M$ is finite (equivalently, bounded).  %for such $k$, each $Q_k$ is a finite set, that is, $Q_k=\{q^k_1,\ldots,q^k_{s(k)}\}$ for some $s(k)\in \bn$. Again by our assumption, 
Then there is the greatest element $\ov{q}$ in $M$ and $p(x_q,x_q)\geq p(x,x)+\varepsilon$ for $q>\ov{q}$. By the choice of $\varepsilon$ and since $\ov{q}\in M$, it holds $p(x_{\ov q},x_{\ov q})< p(x,x)+\varepsilon<\rho_p+\varepsilon_1$. Hence, by condition \pref{cnd:1:2}, for $\varepsilon'\in (0,p(x,x)+\varepsilon-p(x_{\ov q},x_{\ov q}))$, there is $n \in \bn$ for which $p(T^{n}(x_{\ov{q}}),T^{n}(x_{\ov{q}}))\leq p(x_{\ov q},x_{\ov q})+\varepsilon'$. But this implies that $p(x_{n+\ov{q}},x_{n+\ov{q}})=p(T^{n}(x_{\ov{q}}),T^{n}(x_{\ov{q}}))\leq p(x_{\ov q},x_{\ov q})+\varepsilon'< p(x,x)+\varepsilon$, while $n+\ov{q}>\ov{q}$, which is impossible. Therefore, the set $M\subset \bn$ is infinite. 
We conclude that, for each $\varepsilon>0$, there is arbitrarily large positive integer $q_\varepsilon$ for which $p(x_{q_\varepsilon},x_{q_\varepsilon})< p(x,x)+\varepsilon$. But this, in view of the fact that the sequence $(p(x_q,x_q))_{q\in \bn_0}$ converges to $r_x$, implies that $r_x\leq p(x,x)$.
\end{proof}

We are now ready to state and prove the main result of this paper. 

\begin{thrm}\label{thm:1} 
Let $(U,p)$ be a $\rho_p$-complete partial metric space. Suppose that a mapping $T:U\to U$ satisfies condition \pref{cnd:1:1} for some $\alpha\in [0,1)$ and condition \pref{cnd:1:2} for some $\varepsilon_1>0$, then $U_p\neq\emptyset$ and there exists a unique $\ov{x}\in U_p$ for which $T(\ov{x})=\ov{x}$. Moreover, $\lim_{q\to\infty}p(\ov{x},T^q(x))=\lim_{q,q'\to\infty}p(T^q(x),T^{q'}(x))=p(\ov{x},\ov{x})$ for any $x\in U_p$. 
\end{thrm}

\begin{proof}
Let us first prove that $U_p\neq 0$ by showing that there exists a $\rho_p$-Cauchy sequence in $U$. Let $\varepsilon_1>0$ be as in \pref{cnd:1:2}. By Lemma \ref{lm:1}, we may assume that, for any $q\in \bn$ with $1/q<\varepsilon_1$, there exists $y_q\in U$ for which $\rho_p\leq p(T^m(y_q),T^n(y_q))\leq \rho_p+1/q<\rho_p+\varepsilon_1$, $m,n\in \bn_0$. We claim that the sequence $(y_q)_{q\in \bn}$ is $\rho_p$-Cauchy. Indeed, let us fix $\varepsilon\in (0,\varepsilon_1)$ and observe that, for any $q,q'\in \bn$ with $1/q<\varepsilon,\, 1/{q'}<\varepsilon$, condition \pref{cnd:1:2} implies the existence of $n_\varepsilon\in \bn $ for which $p(T^{n_\varepsilon}(y_q),\,T^{n_\varepsilon}(y_{q'}))<\max\{p(y_q,y_q),p(y_{q'},y_{q'})\}+\varepsilon\leq \rho_p+2\varepsilon$. Hence, for such $q,q'$ and the corresponding $n_\varepsilon$,
\begin{multline*}
p(y_q,y_{q'})\leq p(y_q,T^{n_\varepsilon}(y_{q}))+p(T^{n_\varepsilon}(y_q),T^{n_\varepsilon}(y_{q'}))+p(T^{n_\varepsilon}(y_{q'}),y_{q'})-p(T^{n_\varepsilon}(y_q),T^{n_\varepsilon}(y_q))\\-p(T^{n_\varepsilon}(y_{q'}),T^{n_\varepsilon}(y_{q'}))\leq 1/q+\max\{p(y_q,y_q),p(y_{q'},y_{q'})\}+\varepsilon+1/{q'}<\rho_p+4\varepsilon,
\end{multline*}
which proves that the sequence is $\rho_p$-Cauchy. Therefore, $U_p\neq \emptyset$.
To prove the existence of a fixed point of $T$ in $U_p$ it suffices to show that the assumptions of Lemma~\ref{lm:2} are satisfied for a point $x\in U_p$. By Lemma \ref{lm:1}, the sequence $(T^q(x))_{q\in \bn}$ is $\rho_p$-Cauchy for any $x\in U_p$. So, suppose that for a given point $x\in U_p$ and $\ov{x}\in U_p$, it holds $\lim_{q\to\infty} p(T^q(x),\ov x)=\lim_{q\to\infty} p(T^q(x),T^q(x))=p(\ov x,\ov x)$. Similarly, since $\ov{x}\in U_p$, there exists $\ov{x}'\in U_p$ for which $\lim_{q\to\infty} p(T^q(\ov x),\ov x')=\lim_{q\to\infty} p(T^q(\ov x),T^q(\ov x))=p(\ov x',\ov x')$. Let us fix $\varepsilon\in (0,\varepsilon_1)$ and by $M$ denote the set of $q\in \bn$ for which $p(T^{q}(x),T^{q}(\ov{x}))<\rho_p+\varepsilon$. For any $\varepsilon'\in (0,\varepsilon)$, for large $q\in\bn$, $p(T^q(x),T^q(x))<\rho_p+\varepsilon'$ and $p(T^q(\ov{x}),T^q(\ov{x}))<\rho_p+\varepsilon'$. Reasoning as in the proof of Lemma \ref{lm:1} (with small modifications due to the possibility that $x\neq \ov{x}$ - use \pref{cnd:1:2} in such a case) we conclude that the set $M\subset \bn$ is infinite. Hence, for any $\varepsilon>0$, there is $n_\varepsilon\in \bn$ with $$p(T^{n_\varepsilon}(x),\ov{x})\leq \rho_p+\varepsilon,\, p(T^{n_\varepsilon}(\ov{x}),\ov{x}')\leq \rho_p+\varepsilon\text{ and }p(T^{n_\varepsilon}(x),T^{n_\varepsilon}(\ov{x}))\leq  \rho_p+\varepsilon.$$ Therefore \begin{multline*}\rho_p\leq p(\ov{x},\ov{x}')\leq p(\ov x,T^{n_\varepsilon}(x))+p(T^{n_\varepsilon}(x),T^{n_\varepsilon}(\ov x))+p(T^{n_\varepsilon}(\ov x),\ov{x}')-p(T^{n_\varepsilon}({x}),T^{n_\varepsilon}({x}))\\-p(T^{n_\varepsilon}(\ov{x}),T^{n_\varepsilon}(\ov{x}))\leq\rho_p+3\varepsilon.\end{multline*} This reveals that $p(\ov{x},\ov{x}')=\rho_p$. Thus, $\ov{x}=\ov x'$ and by Lemma \ref{lm:2}, $\ov{x}=T(\ov{x})$. Uniqueness of the fixed point in $U_p$ easily follows from condition (\ref{cnd:1:2}). By the above facts, we get $\lim_{q\to\infty}p(\ov{x},T^q(x))=\lim_{q,q'\to\infty}p(T^q(x),T^{q'}(x))=p(\ov{x},\ov{x})$ for any $x\in U_p$. The proof is now complete.
\end{proof} 

\section{Examples and remarks} \label{sec:discussion}

Let us first relate Theorem \ref{thm:1} to some results known from the literature. 

\begin{rmk}\label{rmk:1} 
\textup{Theorem \ref{thm:1} generalizes Theorem 3.1 in \cite{IPR}. Indeed, if we suppose, as in \cite{IPR}, $(U,p)$ is a complete partial metric space and $T:U\to U$ is such that \begin{equation}\label{eq:6} p(T(x),T(y))\leq \max\{\alpha p(x,y), p(x,x), p(y,y)\},\end{equation} $x,y\in U$, for a fixed $\alpha\in (0,1)$, then we have, for $x\in U$ and $x_q:=T^q(x)$, $q\in \bn_0$, 
\begin{multline*}p(x_1,x_{q+1})\leq \max\{\alpha p(x_{0},x_{q}),p(x_{0},x_0), p(x_{q},x_{n})\}\\\leq
\max\{\alpha p(x_0,x_{q}),p(x_0,x_0), \alpha p(x_{q-1},x_{q-1}),p(x_{q-1},x_{n-1})\}\\=\max\{\alpha p(x_0,x_{q}),p(x_0,x_0),p(x_{q-1},x_{q-1})\}
\leq\ldots\leq 
\max\{\alpha p(x_0,x_{q}),p(x_0,x_0)\},
\end{multline*}
which shows that condition \pref{cnd:1:1} is satisfied. To see that condition \pref{cnd:1:2} holds it suffices to observe that, for $x,y\in U,\,x_q:=T^q(x),\, y_q:=T^q(y), \, q\in \bn_0$, by \pref{eq:6} we have
\begin{multline*}p(x_{q},y_{q})\leq \max\{\alpha p(x_{q-1},y_{q-1}),p(x_{q-1},x_{q-1}), p(y_{q-1},y_{q-1})\}\\\leq\max\{\alpha^2 p(x_{q-2},y_{q-2}), p(x_{q-2},x_{q-2}), p(y_{q-2},y_{q-2})\}\leq\ldots\\\leq 
\max\{\alpha^q p(x_0,y_0), p(x_0,x_0), p(y_0,y_0))\}.
\end{multline*}}
\end{rmk} 

\begin{rmk}\label{rmk:2} 
\textup{Let us also notice that condition \pref{cnd:1:1} is true for mappings satisfying rather general assumptions of main results of \cite{CKT}. Namely, if $(U,p)$ is a complete partial metric space and $T:U\to U$ meets conditions of Theorem 2.1 in \cite{CKT}, then for some $\alpha\in [0,1)$ we have: 
\begin{multline*}\label{eq:7} p(T(x),T(y))\\\leq \max\{\alpha p(x,y), \alpha p(x,T(x)), \alpha p(y,T(y)), \frac{\alpha}{2}[p(x,T(y))+p(y,T(x))], p(x,x),p(y,y)\},\end{multline*} $x,y\in U$.
Thus, for $x\in U$, if $x_q:=T^q(x),\, q\in \bn_0$, then, for $q\in\bn$,
\begin{equation}\label{eq:4a}p(x_q,x_{q})\leq \max\{\alpha p(x_{q-1},x_{q}), p(x_{q-1},x_{q-1})\},\end{equation}
and
\begin{multline}\label{eq:4b}p(x_q,x_{q+1})\leq \max\{\alpha p(x_{q-1},x_{q}), \alpha p(x_{q},x_{q+1}), p(x_{q-1},x_{q-1}), p(x_{q},x_{q})\}\\\leq\max\{\alpha p(x_{q-1},x_{q}), \alpha p(x_{q},x_{q+1}), p(x_{q-1},x_{q-1})\}\leq \max\{\alpha p(x_{q-1},x_{q}), p(x_{q-1},x_{q-1})\}\\\leq p(x_q,x_{q-1}),\end{multline}
which result in $p(x_q,x_{q+1})\leq p(x_{q-1},x_{q})\leq\ldots\leq p(x_1,x_0)$. This implies $p(x_q,x_q)\leq \max\{\alpha p(x_{q-1},x_q),\linebreak p(x_{q-1},x_{q-1})\}\leq \max\{\alpha p(x_{q-1},x_q),\alpha p(x_{q-2},x_{q-1}), p(x_{q-2},x_{q-2})\}\leq \ldots \leq \max\{\alpha p(x_0,x_1),p(x_0,x_0)\}.$ It is clear that that condition (\ref{cnd:1:1}) is satisfied for $q=1$. Suppose that it works for some $q\geq 1$. But that it also holds for $q+1$ flows easily from the following inequalities
\begin{multline*}p(x_1,x_{q+2})\\\leq \max\{\alpha p(x_0,x_{q+1}), \alpha p(x_0,x_1), \alpha p(x_{q+1},x_{q+2}),\frac{\alpha}{2}[p(x_0,x_{q+2})+p(x_{q+1},x_1)], p(x_0,x_0),p(x_{q},x_{q})\}\\\leq \max\{\alpha p(x_0,x_{q+2}), \alpha p(x_0,x_{q+1}), \alpha p(x_0,x_1),p(x_0,x_0),\alpha p(x_{q+1},x_1)\}
\end{multline*}
and the inductive assumption.
It follows that condition (\ref{cnd:1:1}) is satisfied.}
\end{rmk}

Now, we are going to consider a few examples which, roughly speaking, illustrate character of conditions \pref{cnd:1:1} and \pref{cnd:1:2}. 
\begin{xmpl}\label{xmpl:1} 
\textup{Let, for $x,y\in U:=\br$, $$p(x,y):=\left\{\begin{array}{cc}|x-y|,& x,y\leq 0,\\|x-y|+1,& \text{otherwise}.\end{array}\right.$$
The pair $(U,p)$ is a $0$-complete partial metric space and $U_p=\{x\in \br:\, x\leq 0\}$. Let us define a mapping $T:U\to U$ by 
$$T(x):=\left\{\begin{array}{cc}\frac{1}{2}x,& x\leq 0,\\ 1,& x>0.\end{array}\right.$$
Since $T^q(x)=\frac{1}{2^q}x,\, x\leq 0$, we have $p(T(x),T^{q+1}(x))=p(2^{-1}x,2^{-(q+1)}x)=2^{-1}|x|(1-2^{-q})\leq\frac{3}{4} p(x,T^{q}(x))=\frac{3}{4}p(x,2^{-q}x)=\frac{3}{4}|x|(1-2^{-q})$, $x\leq 0$, $q\in\bn$. Further, if $x>0$, then $p(T(x),T^{q+1}(x))=p(1,1)=1\leq \max\{\frac{3}{4}p(x,T^{q}(x)),p(x,x)\}=\max\{\frac{3}{4}(1+|x-1|),1\}$, $q\in \bn$. Hence, mapping $T$ satisfies condition \pref{cnd:1:1} with $\alpha=\frac{3}{4}$. Condition \pref{cnd:1:2} is fulfilled for $\varepsilon_1=\frac{1}{2}$. The claim of Theorem \ref{thm:1} holds true and the unique fixed point in $U_p$ is $x=0$. Moreover, $T(1)=1$, but $1\notin U_p$. \newline
Let $M(x,y):=\max\{p(x,y),p(x,T(x)),p(y,T(y)),p(x,T(y)),p(y,T(x))\},\, x,y\in U$. Observe that for, $y\in(0,1)$,  $p(T(0),T(y))=p(0,1)=2>1$, $p(0,0)=0$, $p(y,y)=1$, and
\begin{multline*}M(0,y)=\max\{p(0,y),p(0,T(0)),p(y,T(y)),p(0,T(y)),p(y,T(0))\}\\=\max\left\{p(0,y),p(0,0),p(y,1),p(0,1),p(y,0)\right\}=\max\{1+y,0, 2-y,2,1+y\}=2.
\end{multline*} 
Therefore, there is no $\alpha\in (0,1]$ or a function $\varphi:[0,\infty)\to[0,\infty)$, $\varphi(t)<t,\, t>0$, for which $p(T(x),T(y))\leq \max\{\alpha M(x,y),p(x,x),p(y,y)\}$ for $x,y\in U$, or $$p(T(x),T(y))\leq \max\{\varphi(M(x,y)),p(x,x),p(y,y)\},$$ $x,y\in U$. This implies that results of \cite{CKT}, \cite{IPR}, \cite{R} and \cite{ZZ} are not applicable for the mapping $T$, even if we weaken the requirement of completeness to $\rho_p$-completeness in the assumptions of theorems in those papers.}

\textup{
Similar observation refers to the main results of \cite{IPR3}. Namely, for $q\in \bn$, we have $T^q(0)=0$, and $T^{q}(\frac{1}{2})=1$. Consequently, for $j\in \{1,\ldots,q\}$, $p(0,T^q(\frac{1}{2}))=p(0,1)=2=p(T^q(0),T^q(\frac{1}{2}))$. Therefore, because $p(0,0)=0$, $p(\frac{1}{2},\frac{1}{2})=1$ and $p(0,\frac{1}{2})=\frac{3}{2}$, there is no $\alpha \in [0,1)$ such that for some $q\in \bn$ it would be $p(T^q(0),T^{q}(\frac{1}{2}))\leq  \max\{\alpha p(0,\frac{1}{2}),\, \alpha p(0,T(\frac{1}{2})),\ldots, \alpha p(0,T^q(\frac{1}{2})),p(T^{q-1}(\frac{1}{2}), T^{q-1}(\frac{1}{2})), p(0,0)\}$. Hence, $T$ is neither $C$-operator nor $C_1$-operator, as defined in \cite{IPR3}, and main results of that paper are not applicable for $T$.}
\end{xmpl} 
\begin{xmpl}\label{xmpl:2} 
\textup{Let $(U,p)$ and $T$ be as in Example \ref{xmpl:1}. Define a mapping $T_1:U\to U$ by 
$$T_1(x):=\left\{\begin{array}{cc}T(x),& x\neq 0,\\ 1,& x=0.\end{array}\right.$$
We have, for $q\in \bn$, $p(T_1(0),T_1^{q+1}(0))=p(1,1)=1\leq \frac{3}{4}p(0,T_1(0))=\frac{3}{4}p(0,1)=\frac{3}{2}$, which in view of Example \ref{xmpl:1} and the definition of $T_1$ shows that condition \pref{cnd:1:1} is satisfied with $\alpha=\frac{3}{4}$. Observe that $T_1(0)=1$, so there is no fixed point of $T_1$ in $U_p=(-\infty,0]$, though $T_1(1)=1$. From equalities $T_1^q(0)=1$, $T^q_1(x)=2^{-q}x$, $x<0$, we get $p(T_1^q(0),T_1^q(x))=p(1,2^{-q}x)=2-2^{-q}x>0=\rho_p$, $q\in\bn$, which directly reveals that condition \pref{cnd:1:2} is not satisfied. We could derive the same conclusion from Theorem \ref{thm:1}.}
\end{xmpl}
\begin{xmpl}\label{xmpl:3}\textup{Let, for $x,y\in U:=\{0,1\}$, $$p(x,y):=\left\{\begin{array}{cc}1+|x-y|,& x\neq 0\text{ or } y\neq 0,\\0,& x=y=0.\end{array}\right.$$
The space $(U,p)$ is a $0$-complete partial metric space  and $U_p=\{0\}$. Let us define a mapping $T:U\to U$ by 
$$T(x):=\left\{\begin{array}{cc}0,& x=1,\\ 1,& x=0.\end{array}\right.$$
The mapping $T$ does not satisfy the condition \pref{cnd:1:1} for any $\alpha\in [0,1)$, because $p(T(0),T^{2q}(0))=p(0,T^{2q+1}(0))=p(0,1)=2$ and $p(0,T^{2q}(0))=p(0,0)=0$, $q\in \bn_0$. It is obvious that condition \pref{cnd:1:2} is satisfied with $\varepsilon_1=1/2$. Remark that $T$ has no fixed points while $T^{2q}(0)=0,\,T^{2q}(1)=1$, $q\in \bn_0$. It may be shown that, the mapping $\ov{T}:U\to U$, $\ov T(x):=T^2(x),\, x\in U$, satisfies both conditions \pref{cnd:1:1} and \pref{cnd:1:2} (with $\alpha=1/2,\, \varepsilon_1=1/2$).}
\end{xmpl}
\begin{xmpl}\label{xmpl:4}\textup{Let, for $m,n\in \bn_0$, $$p(m,n):=\left\{\begin{array}{cc}1/m+1/n,& m,n\in \bn,m\neq n\\0,& m=n,\\1/\max\{m,n\},& otherwise.\end{array}\right.$$
The space $(\bn_0,p)$ is a compact metric space and, as such, it is also a $0$-complete partial metric space. Let us define a mapping $T:\bn_0\to\bn_0$ by 
$T(m):=m+1$, $m\in\bn_0$. Notice, that $T^q(m)=q+m,\, m\in\bn_0,\, q\in \bn$. Hence, $p(T^q(m),T^q(m'))\to 0$ for $m,m'\in \bn_0$, as $q\to\infty$. This shows that condition (\ref{cnd:1:2}) is satisfied. It is obvious that $T(m)\neq m$, $m\in \bn$. 
So, even if $(U,p)$ is a compact metric space and $T:U\to U$ meets condition (\ref{cnd:1:2}) there is no guarantee that $T$ has a fixed point. Notice that any iterate $T^q$ of $T$ has no fixed point, though it meets condition (\ref{cnd:1:2}), $q\in \bn$.
}
\end{xmpl} 
\begin{xmpl}\label{xmpl:5}\textup{Let $U:=\bn_0\cup\{n+\frac{1}{2k}:\, n\in\bn_0,\, k\in \bn\}$ and define a metric $p$ on $U$ by $p(x,y):=|x-y|$, $x,y\in U$. It is clear that $(U,p)$ is a complete metric space. We shall now define a mapping $T:U\to U$ satisfying condition (\ref{cnd:1:1}) but not (\ref{cnd:1:2}). We put $$T(x):=\left\{\begin{array}{cc}n+1+\frac{1}{2},& x=n\in \bn_0,\\ n+\frac{1}{2}\frac{1}{2k},& x=n+\frac{1}{2k},\,n\in \bn_0,\, k\in \bn.\end{array}\right.$$
Obviously, $T(x)\neq x,\, x\in U$. 
Observe that, for $x\in \bn_0$, we have $T^{q+1}(x)=x+1+\frac{1}{2^{q+1}}$, $p(T(x),T^{q+1}(x))=\frac{1}{2}-\frac{1}{2^{q+1}}<\frac{1}{2}$ and $p(x,T^q(x))=1+\frac{1}{2^{q}}\geq 1$, $q\in \bn$. Thus, for $x\in \bn_0$, $q\in \bn$, we have $p(T(x), T^{q+1}(x))\leq \frac{3}{4}\max\{p(x,T^{q+1}),\ldots,p(x,T(x))\}$. 
Now, if $x=n+\frac{1}{2k}$, for some $n\in \bn_0,\, k\in \bn$, we get $T^{q+1}(x)=n+\frac{1}{2^{q+1}k}$, $p(x,T^{q+1}(x))=\frac{1}{2k}(1-\frac{1}{2^{q+1}})\geq \frac{1}{2k}\frac{3}{4}$ and $p(T(x),T^{q+1}(x))=\frac{1}{2k}(\frac{1}{2}-\frac{1}{2^{q+1}})\leq \frac{1}{2k}\frac{1}{2}$, $q\in \bn$. Hence we obtain $p(T(x),T^{q+1}(x))\leq \frac{3}{4}p(x,T^{q+1}(x))$, $q\in \bn$. Therefore, for $x=n+\frac{1}{2k}$ and $q\in \bn$ we have $p(T(x), T^{q+1}(x))\leq \frac{3}{4}\max\{p(x,T^{q+1}),\ldots,p(x,T(x))\}$. We conclude that the mapping $T$ satisfies condition (\ref{cnd:1:1}). That condition (\ref{cnd:1:2}) is violated follows from the fact that $p(T^{q}(0), T^q(1))=1$ for $q\in \bn$.
It is also clear that any iteration $T^q$ of $T$ has no fixed point while it meets condition (\ref{cnd:1:1}) with $\alpha=\frac{3}{4}$, $q\in \bn$.
Observe that the example shows that condition (\ref{cnd:1:1}) itself is not enough to ensure the existence of a fixed point even if the space $(U,p)$ is a complete metric space.
}
\end{xmpl} 

\begin{rmk}\label{rmk:3}
\textup{In contrast to Example \ref{xmpl:5} it holds that if $(U,p)$ is a complete metric space and $T:U\to U$ satisfies condition (\ref{cnd:1:1}) and is orbitally continuous (that is, for any $x\in U$, $(T(T^{q_j}(x)))_{j\in \bn}$, where $(q_j)_{j\in \bn}$ is a subsequence of $(q)_{q\in \bn}$, converges to $T(x')$ whenever $(T^{q_j}(x))_{j\in \bn}$ converges to $x'\in U$ as $j\to+\infty$; see \cite{Cir}), then $T$ possesses a fixed point. This stems easily from Lemma \ref{lm:2a} and the fact that every complete metric space is a $0$-complete partial metric space. Nevertheless, the possibility of occurrence of many fixed points, even infinitely many, is not excluded in that case.}
\end{rmk}


\begin{thebibliography}{99}
\bibitem{ASS} I. Altun, F. Sola and H. Simsek, \emph{Generalized contractions on partial metric spaces}, Topology Appl., 157(2010),  2778--2785. 
\bibitem{BMW} D. Bugajewski, P. Ma\'ckowiak and R. Wang, \emph{On compactness and fixed point theorems in partial metric spaces}, Fixed Point Theory 23(2022), 163--178. 
\bibitem{BW} D. Bugajewski, R. Wang, \emph{On the topology of partial metric spaces}, Math. Slovaca, 70(2020), 135--146.
\bibitem{BKMP} M. Bukatin, R. Kopperman, S. Matthews and H. Pajoohesh, \emph{Partial metric spaces}, Amer. Math. Monthly 8(2009), 708-718.
\bibitem{CKT} K.P. Chi, E. Karapinar, T.D. Thanh, \emph{A generalized contraction principle in partial metric spaces}, Math. Comput. Modelling, 55(2012), 1673--1681.
\bibitem{Cir} L. \'Ciri\'c, \emph{On contraction type mappings}, Math. Balkanica, 1(1971), 52--57.
\bibitem{HRS} R.H. Haghi, Sh. Rezapour and N. Shahead, \emph{Be careful on partial metric fixed point results}, Topology Appl., 160(2013), 450--454.
\bibitem{IPR} D. Ili\'{c}, V. Pavlovi\'{c} and V. Rako\v{c}evi\'{c}, \emph{Some new extensions of Banach's contraction principle to partial metric space}, Appl. Math. Lett., 24(2011), 1326--1330.
\bibitem{IPR2} D. Ili\'{c}, V. Pavlovi\'{c} and V. Rako\v{c}evi\'{c}, \emph{Extensions of the Zamfirescu theorem to partial metric spaces}, Math. Comput. Modelling, 55(2012), 801--809.
\bibitem{IPR3} D. Ili\'{c}, V. Pavlovi\'{c} and V. Rako\v{c}evi\'{c}, \emph{Fixed points of mappings with a contractive iterate at a point in partial metric spaces}, Fixed Point Theory Appl. 2013, 335(2013), 1--18.

\bibitem{Mack2023} P. Ma\'ckowiak, \emph{A Converse of the Banach Contraction Principle for Partial Metric Spaces and the Continuum Hypothesis}, Results Math. 79, 41(2024), https://doi.org/10.1007/s00025-023-02072-5
\bibitem{M} S.G. Matthews, \emph{Partial metric topology}, in: Papers on General Topology and Applications, Flushing, NY, 1992, in: Ann. New York Acad. Sci., vol.728, New York Acad.Sci., New York, 1994, 183--197.
\bibitem{R} S. Romaguera, \emph{Fixed point theorems for generalized contractions on partial metric spaces}, Topology Appl., 159(2012), 194--199.

\bibitem{ST} M.B. Smyth and R. Tsaur, \emph{Hyperconvex semi-metric spaces}, Topology Proc., 2006(2001-2002), 791-810.
\bibitem{ZZ} Y. Zhang, Z. Jiang, \emph{New nonlinear contraction principle in partial metric spaces}, J. Funct. Spaces 2016, Article ID 8404529 (2016), 1--9.

\end{thebibliography}
\end{document}